\documentclass[11pt]{article}

\usepackage[dvips]{graphicx}
\usepackage{bm}
\usepackage{amsthm,citesort,algorithm,algorithmic,url}
\newtheorem{thm}{Theorem}
\newtheorem{lemma}{Lemma}
\newtheorem{teigi}{Definition}
\newtheorem{fact}{Fact}

\newtheorem{observation}{Observation}

\floatname{algorithm}{Procedure}

\title{Enumerating Constrained Non-crossing Minimally Rigid Frameworks}
\author{\setcounter{footnote}{0}%
\def\thefootnote{\arabic{footnote}}
David Avis \footnotemark[1]\hspace{3mm}
Naoki Katoh \footnotemark[2]\hspace{3mm}
Makoto Ohsaki \footnotemark[2]\hspace{3mm}
\\
\def\thefootnote{\arabic{footnote}}
Ileana Streinu \footnotemark[3]\hspace{3mm}
Shin-ichi Tanigawa \footnotemark[2]
}

\begin{document}
\maketitle

\begin{abstract}
In this paper we present an algorithm for enumerating 
without repetitions all the non-crossing generically minimally 
rigid bar-and-joint frameworks under edge constraints
(also called constrained non-crossing Laman frameworks) on a given generic set of $n$ points.
Our algorithm is based on the reverse search paradigm of Avis and Fukuda.
It generates each output graph in $O(n^4)$ time and $O(n)$ space, or, slightly different implementation,
in $O(n^3)$ time and $O(n^2)$ space.
In particular, we obtain that the set of all the constrained non-crossing Laman frameworks 
on a given point set is connected by flips
which restore the Laman property.

\vspace{\baselineskip}
\end{abstract}

{
\setcounter{footnote}{0}
\def\thefootnote{\arabic{footnote}}
\footnotetext[1]{School\ of\ Computer Science, McGill University,
Canada.}
\footnotetext[2]{Department of Architecture and Architectural
Engineering, Kyoto University Katsura, Nishikyo-ku, Kyoto 615-8450
Japan, \url{{ohsaki,naoki,is.tanigawa}@archi.kyoto-u.ac.jp}.
Supported by JSPS Grant-in-Aid for Scientific Research on priority areas of New Horizons in Computing.}
\footnotetext[3]{Dept. of Comp. Science, Smith College,
Northampton, MA 01063, USA, \url{streinu@cs.smith.edu}. Supported
by NSF grant CCF-0430990 and NSF-DARPA CARGO CCR-0310661.}
}

\section{Introduction}
\label{introduction}

Let $G$ be a graph with vertices $\{1,\dots,n\}$ and $m$ edges.
$G$ is a {\it minimally rigid graph} (also called {\em Laman graph}) if $m=2n-3$ and every subset of $n'\leq n$ 
vertices spans at most $2n'-3$ edges. 
An embedding of the graph $G(P)$ on a set of points $P=\{ p_1, \cdots, p_n\}\subset R^2$ 
is a mapping of the vertices to points in the Euclidean
plane $i\mapsto p_i\in P$. The edges $ij$ of $G$ are mapped to
straight line segments $p_ip_j$. 
An embedding is {\it planar} or {\it non-crossing} if no pair of segments $p_i p_j$ and
$p_kp_l$ corresponding to non-adjacent edges $i,j\not\in\{k,l\}$ have a point in common.

An embedded Laman graph on the {\em generic} point sets is called {\em Laman framework} and has the
special property of being {\em minimally rigid}
\cite{laman:Rigidity:1970,graver:servatius:servatius:CombinatorialRigidity:1993},
when viewed as bar-and-joint {\em frameworks} with fixed edge-lengths,
which motivates the tremendous interest in their properties.
Let $F$ be a set of non-crossing edge (bar) on $P$. 
A Laman framework containing $F$ is called {\em F-constrained}.
In this paper we give an algorithm for enumerating all the {\em
F-constrained non-crossing Laman frameworks}.

In this paper we only consider planar embedded Laman graphs, namely Laman frameworks.
Then we simply denote a vertex $p_i$ by $i$ and an edge $p_ip_j$ by $ij$.

\medskip
\noindent {\bf Novelty.} \hspace{1em} To the best of our knowledge, this is the
first algorithm proposed for enumerating (without repetitions, in
polynomial time and without using additional space) all the $F$-constrained non-crossing
Laman frameworks. We achieve $O(n^3)$ time per graph in $O(n^2)$ space by using reverse search.
For the unconstrained case, using a different method,
the same time and space complexity results we obtained in \cite{avis:katoh:enumerate-Laman:2006}.
The method presented here uses matroid property and is considerably simpler than \cite{avis:katoh:enumerate-Laman:2006}.

\medskip
\noindent {\bf Historical Perspective.}  \hspace{1em} The reverse search
enumeration technique of Avis and Fukuda
\cite{avis:fukuda:vertex-enumeration:1992,avis:fukuda:reverse-search:1996}
has been successfully applied to a variety of combinatorial and
geometric enumeration problems.
The necessary ingredients to use the
method are an implicitly described connected graph
on the objects to be generated, and an implicitly defined
spanning tree in this graph.
In this paper we supply these
ingredients for the problem of generating Laman frameworks.

Relevant to the historical context of our work are the results of
Bereg \cite{bereg:enumerating-triangulations:2002,bereg:enumerating-pseudo-triangulations:2005}
using reverse search combined with data-specific lexicographic
orderings to enumerate triangulations and pointed
pseudo-triangulations of a given point set. We notice in passing
that there exist several other algorithms for enumerating
(pseudo-)triangulations
\cite{dumitrescu:gaertner:pedroni:welzl:enumeratingTriangPaths:2001,bronnimann:kettner:pocchiola:snoeyink:enumeratingPointedPseudoTriang:2005,streinu:aichholzer:rote:speckmann:wads2003},
but they are based on different techniques.

Also relevant is the pebble game algorithm of Jacobs and
Hendrickson \cite{jacobs:hendrickson} for 2-dimensional rigidity,
see also \cite{berg:jordan:2003}. Our complexity analysis relies
on the recent results, due to Lee, Streinu and Theran
\cite{lee:streinu:2005,lee:streinu:theran:2005}, regarding the
detailed data-structure complexity of finding and maintaining
rigid components during the pebble game algorithm. Indeed, the
time-space trade-off of our algorithm is inherited from
\cite{lee:streinu:theran:2005}.

\medskip
\noindent {\bf Related previous works.}  \hspace{1em}
To better put our problem (and solution) in context and relate it
to previous work, we briefly discuss now the difference between
generating arbitrary Laman frameworks, as opposed to just
pointed pseudo-triangulations.

A {\it pointed pseudo-triangulation} is a special case of a non-crossing Laman framework \cite{streinu:pseudoT:2005},
where every vertex in the embedding is incident to an angle larger
than $\pi$.
Pseudo-triangulations are connected via simple {\em flips}, in
which the removal of any non-convex-hull edge leads to the choice
of a {\em unique} other edge that can replace it, in order to
restore the pointed pseudo-triangulation property.
The flip graph of all pointed pseudo-triangulations is a connected
subgraph of the graph of all the Laman frameworks. In fact, it is the
one-skeleton of a polytope
\cite{streinu:rote:santos:expansive:2003}, and the reverse search
technique can be directly applied to it. 
Bereg's efficient
algorithm makes use of specific properties of pointed
pseudo-triangulations which do not extend to arbitrary non-crossing Laman frameworks. 
In particular, remove-add flips are {\em not unique}, relative to the removed edge, 
in the case of a non-crossing Laman framework. 

In our previous paper \cite{avis:katoh:enumerate-Laman:2006}
we showed that a set of all non-crossing Laman frameworks on a point set is connected via these flips.
The proof relies on some properties of one-degree-of-freedom mechanisms 
which is obtained when removing one edge from the Laman frameworks,
which cannot extend to edge constrained case.

In this paper we use triangulation as an important tool.
We use the fact that (constrained) non-crossing Laman frameworks form the bases of a matroid
defined on any triangulation of $P$ to define a parent function.
The edge constrained case appears in an application to structural engineering as will be seen below,
and the proposed algorithm will much reduce the output size compared with the unconstrained case.
Therefore, we can enumerate non-crossing Laman frameworks for previously unsolved real problems.


\medskip
\noindent {\bf Motivating application.} \hspace{1em} 
We describe now briefly how this problem came to our attention via the work of the third author.
Graph theoretical approaches are widely used in {\em structural mechanics} \cite{kaveh:structuralMechanics:2004}, 
where the edges and vertices in the graph represent the bars and rotation-free joints of a structure called a {\em truss}.
It is well-known that the stiffest truss under static loads is statically determinate 
that is equivalent to a Laman graph \cite{bendsoe:sigmund:topologyOptimization:2003}.  

Another bar-joint system, which is widely used in industrial application, 
is a link mechanism that is unstable and generates large deformation or changes the direction of the nodal displacement.
The applications of link mechanisms can be found in, e.g., suspension of automobiles, robot hands, umbrellas, crank shafts, etc.
Kawamoto {et al.}\cite{kawamoto:bendsoe:sigmund:planeMechanismDesign:2004}
presented a method by enumeration of the planar graphs to find an optimal mechanism.
However, their method is developed for their specific problem, and no general approach was given.

Recently, a new type of mechanism called compliant mechanism has been developed and applied mainly in the field of micro-mechanics.
A compliant mechanism has flexible parts to stabilize the structure, which is contrary to the conventional unstable mechanism.
Although a compliant mechanism is usually modelled as a continuum with elastic joints,
it is possible to generate the similar mechanism by a bar-joint system.
Ohsaki and Nishiwaki \cite{ohsaki:hishiwaki:compliantMechanisms:2005} presented a method for generating flexible bar-joint mechanisms
using nonlinear programming approach, and found that the optimal
structure is statically determinate, i.e. minimally rigid.
They utilized snapthrough behaviour to generate multi-stable mechanism that has multiple self-equilibrium states.
Such mechanism can be used for switching device, robot hand, gripper, deployable structure, etc.
In their method, the optimal locations of bars and joints are found from the highly connected initial structure
that has bars between all the pairs of the nodes whose distances are small enough.
Combining an implementation of our our earlier method for generating unconstrained
Laman frameworks \cite{avis:katoh:enumerate-Laman:2006}
with this nonlinear programming approach, we could obtain
many new compliant mechanisms with up to 10 joints \cite{katoh:CJK}. However the number of
Laman frameworks grows too rapidly to allow a complete enumeration for
significantly larger examples.

In view of practical requirements, the optimal structure should
not have intersecting members, and some pre-selected members should always exist.
Therefore, the computational cost may be much reduced
if the candidate set of statically determinate non-crossing trusses (non-crossing Laman frameworks) are first enumerated.
Thus, it is desirable to enumerate Laman frameworks for which the existing edges are specified,
and this is the problem considered in the current paper.


\section{Preliminaries}
\label{sec:preliminaries}

Let $L$ be a non-crossing Laman framework on a given point
set $P$. 
%
A {\em mechanism} is a flexible framework obtained by removing one
or more edges from a Laman framework. Its {\it number of
degrees of freedom} or {\em dof}'s, is the number of removed
edges.  We will encounter mostly {\em one-degree-of-freedom (1dof)
mechanisms}, which arise from a Laman framework by the removal of one
edge.
In particular, a mechanism with $k$ dofs has exactly $2n-3-k$
edges, and each subset of $n'$ vertices spans at most $2n'-3$
edges. A subset of some $n'$ vertices spanning {\em exactly}
$2n'-3$ edges is called a {\em rigid block}. A {\em maximal block}
is called a {\em rigid component} or a {\em body}. 

The Laman frameworks on the generic point set with $n$ points form the set of {\em bases} of
the {\em generic rigidity matroid} (or simply called {\em Laman matroid}) on $K_n$, see
\cite{graver:servatius:servatius:CombinatorialRigidity:1993}. 
The bases have
all the same size $2n-3$. Bases may be related via the {\em base
exchange} operation, which we will call a {\em flip} between two
Laman frameworks. Two Laman frameworks $L_1$ and $L_2$ are connected by a
flip if their edge sets agree on $2n-4$ elements. The flip is
given by the pair of edges $(e_1,e_2)$ not common to the two
bases, $e_1\in L_1-L_2$, $e_2\in L_2-L_1$. Using
flips, we can define a graph whose nodes are {\em all} the Laman
frameworks on $n$ vertices, and whose edges correspond to flips. 

It is well-known that the graph whose nodes are the bases of a
matroid connected via flips, is connected. But a priori, the
subset of {\em (constrained) non-crossing Laman frameworks} may not necessarily be. We
will prove this later in Section \ref{sec:main}.
%
%

\medskip
Reverse search is a memory efficient method for visiting all the nodes of
a connected graph that can be defined implicitly by
an adjacency oracle. It can be used whenever a spanning tree
of the graph can be defined implicitly by a {\em parent}
function. This function is defined for each vertex of the
graph except a pre-specified {\it root}.
Iterating the parent function leads to a path to the root from
any other vertex in the graph. The set of such paths defines
a spanning tree, known as the {\it search tree}.

\section{Constrained Non-crossing Laman Frameworks}
\label{sec:main}
Let $T$ be a triangulation on a given set of $n$ points $P$ in the plane 
containing
$k$ triangles. The angle vector of $T$ is the vector of $3k$ interior angles
sorted into non-decreasing order.
Let $F$ be a non-crossing edge set on $P$. 
An {\em F-constrained triangulation} on $P$ is one that contains $F$ as a subset. 
Many facts about these triangulations are contained in the survey by
Bern and Eppstein \cite{BE92}.
If $F$ is an
independent set in the Laman matroid on $K_n$, then 
a Laman framework containing $F$ is called {\em $F$-constrained}.
The following lemma clearly holds from the known fact about matroids (see, e.g.\cite{Welsh}).
\begin{lemma}\label{lemma:minor}
Let $F$ be a non-crossing edge set on $P$ that is an independent set in the Laman matroid on $K_n$.
Let us fix $F$-constrained triangulation $T$ on $P$.
Then a set of bases of the Laman matroid on $T$ whose edge set contains $F$ forms matroid.
\end{lemma}

\begin{lemma}\label{lemma:trilam}
Let $F$ be a non-crossing edge set on $P$ that is an independent set in the 
Laman matroid on $K_n$. Every $F$-constrained
triangulation $T$ on $P$ contains an $F$-constrained 
Laman framework.
\end{lemma}

\begin{proof}
Since $T$ is statically rigid, it is generically rigid,
and hence contains a base $B$ of the Laman matroid (see, e.g.
Whiteley \cite{Whitley97}).
Since $F$ is independent in the Laman matroid, it can be
extended to a base by adding edges from $B-F$.
\end{proof}

Two points $a$ and $b$ are {\em visible} if
no edge of $F$ properly intersects 
the segment $ab$. $ab$ is {\em visible} to
point $c$ if the triangle $abc$ is not properly intersected by an edge of $F$.

\begin{teigi}
(CDT)
An F-constrained Delaunay Triangulation (CDT)
contains the edge $ab$ between points $a$ and $b$ in $P$ if and only
if $a$ is visible to $b$, and some circle through $a$ and $b$ contains
no point of $P$ visible to segment $ab$. We call
$ab$ a D(elaunay)-edge. (Definition 1, \cite{BE92})
\end{teigi}

If $P$ has four or more co-circular points, using a linear transformation
as described in \cite{beichl:sullivan:1995}, we may transform $P$ into a point set $\bar{P}$ with a unique CDT.
$P$ and $\bar{P}$ have the same non-crossing Laman frameworks.
Then we will assume in what follows that $P$ has a unique CDT.

\begin{teigi}
(D-flip)
Let $ac$ be an edge of $T-F$ which is the diagonal of a convex quadrilateral
$abcd$ contained in $T$. 
The edge flip which replaces $ac$ by edge $bd$
is a D-flip if the circumcircle of the triangle $abc$ contains the point
$d$ in its interior. (Equivalently the circumcircle of the triangle $acd$ 
contains the point
$b$ in its interior.) 
We call $ac$ a {\it $F$-illegal edge}, otherwise it is called {\it $F$-legal edge}.
\end{teigi}

Note that a D-flip increases the angle vector lexicographically.
This can be used to prove the following.

 \begin{fact}\label{fact:angle}
The CDT has the lexicographically maximum angle vector of all
$F$-constrained triangulations on $P$. (Theorem 1, \cite{BE92}).
 \end{fact}


\begin{fact}\label{fact:d-flip}
   An $F$-constrained triangulation $T$ can be converted to 
the CDT by at most $O(n^2)$
D-flips, taken in any order. (Lemma 4, \cite{BE92}).
 \end{fact}

Now let us consider non-crossing Laman frameworks.
\begin{teigi}
A \textit{Constrained Delaunay Laman Framework (CDLF)}
is an $F$-constrained Laman framework that is a subset of the $CDT$.
\end{teigi}

\begin{teigi}
(L-flip)
An L-flip is an edge insertion and deletion that takes a Laman framework $L$
to a new Laman framework $L'$.
\end{teigi}

\begin{thm}
\label{theo1}
Every $F$-constrained non-crossing Laman framework $L$ can be transformed to a CDLF by at most $O(n^2)$ L-flips.
\end{thm}
Before giving the proof of Theorem \ref{theo1}, we define the underlying triangulation $T(L)$ for a non-crossing Laman framework $L$,
which is $L$-constrained Delaunay triangulation.
\begin{teigi}
\label{teigi:L-constrained_Delaunay}
Let $L$ be a non-crossing Laman framework.
The $L$-constrained Delaunay triangulation $T(L)$ is constructed by adding edges to $L$
in the following way. 
First, we add convex hull edges if missing in $L$, and
for each planar face of $L$ we compute its (internal)
constrained Delaunay triangulation, and add these new edges.
\end{teigi}
By this construction, we have the following fact:
\begin{fact}
\label{fact:under_tri}
Let $L$ be a $F$-constrained non-crossing Laman framework. 
Then all $F$-illegal edges in $T(L)$ are the edges of $L-F$.
\end{fact}

\begin{proof}[Proof of Theorem 1]
Suppose that $L$ has an $F$-illegal edge.
Then, from Fact \ref{fact:under_tri}, 
such a $F$-illegal edge, say $ac$, is among the edges of $L-F$.
Consider removing $ac$.
The updated underlying triangulation $T(L-ac)$ contains $L-ac$, and 
by Lemma \ref{lemma:trilam} and the matroid property shown by Lemma \ref{lemma:minor}, there is
some edge $st$ in $T(L-ac)$ such that $L'=L-ac+st$ is a non-crossing Laman framework.
The fact that $ac$ is $F$-illegal edge implies that there exists at least one D-flip 
when updating $T(L)$ to $T(L-ac)$,
which lexicographically increases the angle vector.
From Fact \ref{fact:angle} and \ref{fact:d-flip}, repeating this procedure $O(n^2)$ times, we
eventually reach the $F$-constrained Delaunay triangulation $T(L")$, and $L"$ is the
required CDLF.
\end{proof}

For edges $e=ij$ with $i<j$ and $e'=kl$ with $k<l$, 
we use the notation $e\prec e'$ or $e'\succ e$ when $e$ is lexicographically smaller than $e'$ i.e., either $i<k$ or $i=k$ and $j<l$,
and $e=e'$ when they coincide.
For an edge set $A$ we use the notations $\max\{e\mid e\in A\}$ and $\min\{e\mid e\in A\}$ 
to denote the lexicographically largest and smallest labelled edges in $A$, respectively.
\begin{teigi}
(Lexicographic edge list)
Let $E=\{e_1\prec e_2\prec\dots\prec e_m\}$ and $E'=\{e'_1\prec e'_2\prec\dots\prec e'_m\}$ 
be the lexicographically ordered edge lists.
Then $E$ is lexicographically smaller than $E'$ if $e_i\prec e'_i$ for the smallest $i$ such that $e_i\neq e'_i$.
\end{teigi}

\begin{thm}
\label{theo2}
The set of $F$-constrained non-crossing Laman Frameworks on a point set $P$ is connected by 
$O(n^2)$ edge flip operations.
\end{thm}
\begin{proof}
By the previous theorem, from any non-crossing Laman framework we can
perform L-flips $O(n^2)$ time to reach a CDLF, say $L$. Let $L^*$ be the CDLF with
lexicographically smallest edge list.
We show that we can
do edge flips from $L$ to $L^*$, at most $n-3$ times, maintaining the non-crossing Laman property.
Indeed, delete from $L$ the largest indexed edge $ac$ in $L-L^*$. 
By the matroid properties, there will always be an edge in $L^*-L$ to insert that
maintains the Laman framework property. 
In fact, the 1dof mechanism $L-ac$ is not maximum component and can be extended to a base
by adding an edge from $L^*-L$.
Planarity is maintained since
both $L$ and $L^*$ are subsets of the CDT which is non-crossing.

Let $n_e$ be the number of edges of the triangulation on $P$.
By Euler's formula, $n_e=3n-h-3\leq 3n-6$ holds, 
where $h\geq 3$ is the number of points on the convex hull of $P$.
Then $L$ has at most $n-3$ edges which are not in $L^*$.
Since we can replace one of such edges by an edge of $L^*$ by one edge flip,
we eventually reach $L^*$ by $n-3$ edge flips.
\end{proof}

\section{Algorithm}
\label{sec:algorithm}
\begin{figure}[t]
\centering
\includegraphics[scale=0.68]{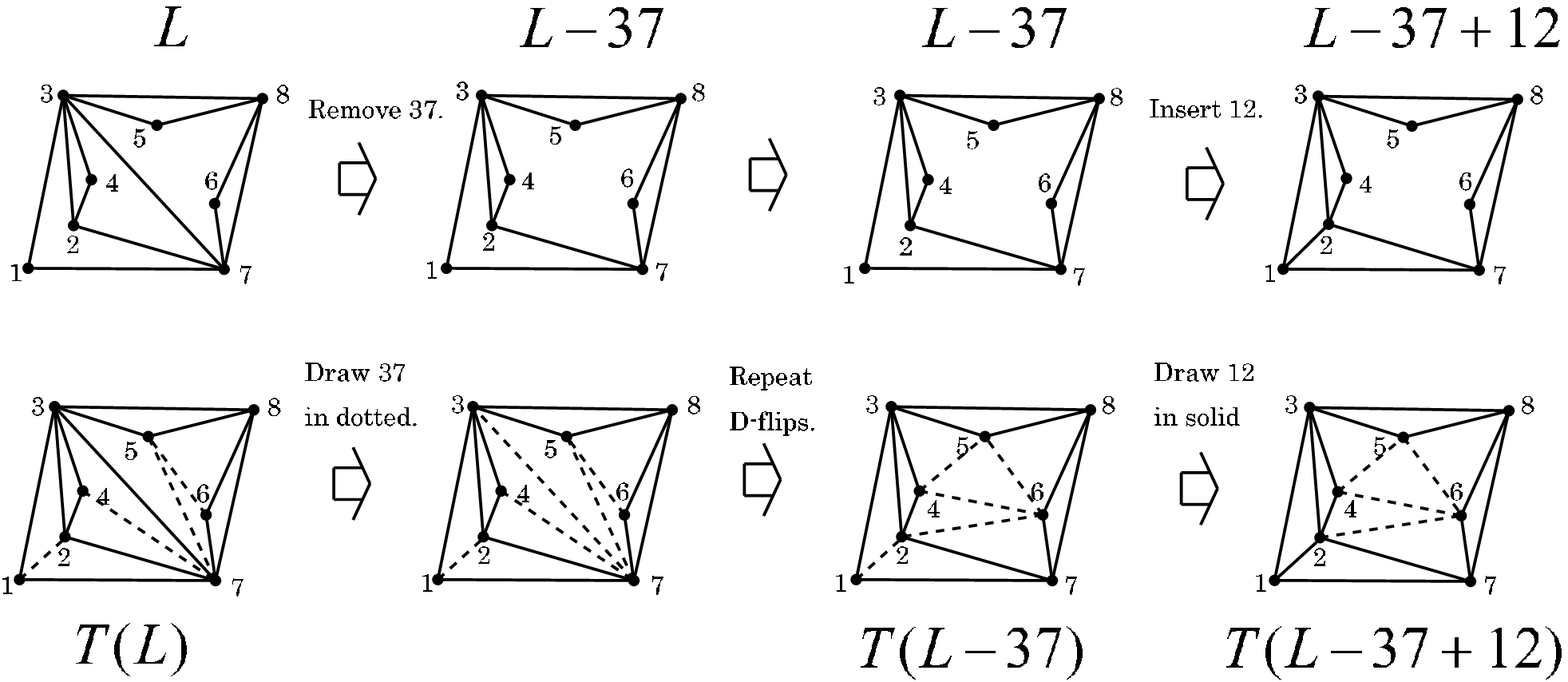}
\caption{The example of the parent function. 
The solid lines represent edges of $L$, and the dotted lines represent additional edges for the underlying triangulations.}
\label{fig:parent}
\end{figure}

Let $\mathcal{L}$ be a set of $F$-constrained non-crossing Laman frameworks on $P$, 
and $\mathcal{DL}$ be a set of the CDLFs.
Clearly $\mathcal{DL}\subseteq \mathcal{L}$ holds.
Let $L^*$ be the CDLF with
lexicographically smallest edge list.
We define the following {\it parent function} $f:\mathcal{L}-\{L^*\}\to\mathcal{L}$ based on Theorem \ref{theo1} and \ref{theo2}.
\begin{teigi}
\label{def:parent}
(Parent function)
Let $L\in \mathcal{L}$ with $L \neq L^*$. 
$L' = L - ac + st$ is the parent of $L$, where \\
Case 1: $L\in\mathcal{DL}$, \\
$ac=\max\{e\mid e\in L-L^*\}$ and $st=\min\{e\in L^*-L\mid L-ac+e\in \mathcal{L}\}$, \\
Case 2: $L\in\mathcal{L}-\mathcal{DL}$,\\
$ac=\max\{e\in T(L)\mid e \mbox{ is $F$-illegal edge in } T(L)\}$ and $st=\min\{e\in T(L-ac)-(L-ac)\mid L-ac+e \in\mathcal{L}\}$.
\end{teigi}
To simplify the notations, we denote the parent function depending on Case 1 and Case 2 
by $f_1:\mathcal{DL}-\{L^*\}\to\mathcal{DL}$ 
and $f_2:\mathcal{L}-\mathcal{DL}\to\mathcal{L}$,
respectively. 
In Fig.\ref{fig:parent}, we show the example of $f_2(L)$ in which $L$ is not a CDLF:
removing the largest indexed $F$-illegal edge $37$, and inserting the
smallest indexed edge $12$ in $T(L-ac)-(L-ac)$,
we obtain another non-crossing Laman framework shown in the rightmost and upper corner of Fig. \ref{fig:parent}.
%
%

\vspace{\baselineskip}

The reverse search algorithm can be considered on the underlying graph in which each vertex corresponds to a non-crossing Laman framework
and two frameworks are {\it adjacent} if and only if one can be obtained from the other by a L-flip.
Then, for $L'\in\mathcal{L}$ the local search is given by an {\em adjacency function}, $Adj$, defined as follows:
\begin{displaymath} 
Adj(L',e_1,e_2):=\left\{ 
\begin{array}{ll} 
L'-e_1+e_2 & \quad \mbox{if }L'-e_1+e_2\in\mathcal{L}, \\ 
null & \quad \mbox{otherwise,}
\end{array} \right. 
\end{displaymath}
where $e_1\in L'-F$ and $e_2\in K_n-L'$.
The number of candidate edge pairs $(e_1,e_2)$ is $O(n^3)$.

Let ${\sf elist}_{L'}$ and ${\sf elist}_{K_n}$ be the list of edges of $L'$ and $K_n$ ordered lexicographically,
let $\delta(L')$ and $\delta(K_n)$ be the number of elements of ${\sf elist}_{L'}$ and ${\sf elist}_{K_n}$
and let ${\sf elist}_{L'}(i)$ and ${\sf elist}_{K_n}(i)$ be the $i$-th elements of ${\sf elist}_{L'}$ and ${\sf elist}_{K_n}$, respectively.
We also denote the above defined adjacency function by $Adj(L',i,j)$ 
for which $e_1={\sf elist}_{L'}(i)$ with $e_1\notin F$ and $e_2={\sf elist}_{K_n}(j)$ with $e_2\notin L'$.
Then, based on the algorithm in \cite{avis:fukuda:vertex-enumeration:1992,avis:fukuda:reverse-search:1996},
we describe our algorithm in Fig. \ref{algo:reverse}.
\begin{figure}[t]
\hrule height 0.8pt
\smallskip
{\bf Algorithm} Enumerating $F$-constrained non-crossing Laman frameworks.
\smallskip
\hrule
\smallskip
\begin{algorithmic}[1]
\STATE $L^*:=${ the CDLF with lexicographically smallest edge list};
\STATE $L':=L^*$; $i,j:=0$; $Output(L')$;
\REPEAT
\WHILE{ $i\leq \delta(L')$}
\STATE $i:=i+1$;
\WHILE{ $j\leq \delta(K_n)$}
\STATE $j:=j+1$;
\IF{${\sf elist}_{L'}(i)\notin F$, ${\sf elist}_{K_n}(j)\notin L'$ and $Adj(L',i,j)\neq$ {\em null}}
\STATE $L:=Adj(L',i,j)$;
\IF{ $f_1(L)=L'$ or $f_2(L)=L'$}
\STATE $L':=L$; $i,j:=0$;
\STATE $Output(L')$;
\STATE {\bf go to} line 4;
\ENDIF
\ENDIF
\ENDWHILE
\ENDWHILE
\IF{$L'\neq L^*$}
\STATE $L:=L'$;
\STATE {\bf if} $L\in {\cal DL}$ {\bf then} $L':=f_1(L)$;  
\STATE {\bf else} $L':=f_2(L)$; 
\STATE determine integers pair $(i,j)$ such that $Adj(L',i,j)=L$;
\STATE $i:=i-1$;
\ENDIF
\UNTIL{$L'=L^*$, $i=\delta(L')$ and $j=\delta(K_n)$};
\end{algorithmic}
\hrule
\caption{Algorithm for enumerating $F$-constrained non-crossing Laman frameworks.}
\label{algo:reverse}
\end{figure}

As we will show later, both the parent function and the adjacency function need $O(n^2)$ time for each process.
Then, the while-loop from line 4 to 17 has $\delta(L')\cdot \delta(K_n)$ iterations
which require $O(n^5)$ time if simply checking the line 8 and 10.
In order to improve $O(n^5)$ time to $O(n^3)$ time we claim the following two lemmas:
\begin{lemma}
\label{lemma:parent1}
Let $L$ and $L'$ be two distinct CDLFs which are subgraphs of $T(L^*)$
for which $L=Adj(L',e_1,e_2)$ for $e_1\in L'-F$ and $e_2\in K_n-L'$. 
Then, $f_1(L)=L'$ holds if and only if $e_1$ and $e_2$ satisfy the following conditions:
\begin{itemize}
\item[(a)] $e_1\in L^{*}$,
\item[(b)] $e_2\in T(L^{*})-L^*$,
\item[(c)] $e_1\prec\min\{e\in L^*-L' \mid L'-e_1+e\in\mathcal{L}\}$,
\item[(d)] $e_2\succ\max\{e\mid e\in L'-L^*\}$.
\end{itemize}
\end{lemma}

\begin{lemma}
\label{lemma:parent2}
Let $L$ and $L'$ be two distinct non-crossing Laman frameworks
for which $L=Adj(L',e_1,e_2)$ for edges $e_1\in L'-F$ and $e_2\in K_n-L'$, and
let $L\in \mathcal{L}-\mathcal{DL}$. 
Then, $f_2(L)=L'$ holds if and only if $e_1$ and $e_2$ satisfy the following conditions:
\begin{itemize}
\item[(a)] $e_1$ is $F$-legal edge in $T(L')$,
\item[(b)] $e_2\in K_n -T(L')$,
\item[(c)] $e_1\prec\min\{e\in T(L')-L'\mid L'-e_1+e\in \mathcal{L}\}$.
\item[(d)] $e_2=\max\{e\in T(L'-e_1+e_2)\mid e \mbox{ is $F$-illegal edge in } T(L'-e_1+e_2)\}$.
\end{itemize}
\end{lemma}

We will explain later (in the proof of Theorem \ref{theo:time}) how Lemmas \ref{lemma:parent1} and \ref{lemma:parent2} are used
to obtain $O(n^3)$ time for generating each output of our algorithm.
Notice that for $L'$ and $L\in \mathcal{L}$ such that $L=L'-e_1+e_2$,
at most one of $f_1(L)=L'$ and $f_2(L)=L'$ holds
from the conditions (b) of Lemmas \ref{lemma:parent1} and \ref{lemma:parent2}.

\begin{proof}[Proof of Lemma \ref{lemma:parent1}]
(``only if''-part.)
Since $f_1(L)=L'$ holds, $e_1$ and $e_2$ must be chosen as $st$ and $ac$ in Case 1 of Definition \ref{def:parent}.
From Definition \ref{def:parent}, $e_2(=ac)\in L-L^*$ holds.
Since $L\in \mathcal{DL}$, $L\subset T(L^*)$ and
$e_2\in T(L^*)-L^*$ holds, thus we have (b).
Similarly since $e_1(=st)\in L^*-L\subset L^*$, we have (a).

From $e_1=st$, we have
\begin{equation}
\label{eq:lem1_eq}
L'-e_1=(L-ac+st)-e_1=L-ac.
\end{equation}
Let $e'=\min\{e\in L^*-L'\mid L'-e_1+e\in \mathcal{L}\}$.
Suppose (c) does not hold. Then $e'\prec e_1$ holds. 
(Note that the equality does not hold since $e_1\in L'-F$.)
Then from Eq.(\ref{eq:lem1_eq}) and $e'\prec e_1=st\prec ac$ (which comes from Definition \ref{def:parent}),
\[
\begin{array}{rcll}
e'&=&\min\{e\in L^*-L'\mid L'-e_1+e\in \mathcal{L}\} & \\
  &=&\min\{e\in L^*-(L-ac+st)\mid L-ac+e\in \mathcal{L}\} & \hspace{1em} \mbox{ (from (\ref{eq:lem1_eq}))} \\
  &=&\min\{e\in L^*-L\mid L-ac+e\in \mathcal{L}\} & \hspace{1em} \mbox{ (from } e'\prec st\prec ac).
\end{array}
\]
Thus, $e'$ would have been selected instead of $e_1$ when the parent function $f_1$ is applied to $L$,
which contradicts $e_1=st$.
Hence, (c) holds.

Finally, let $e''=\max\{e\mid e\in L'-L^*\}$, and suppose that (d) does not hold and $e_2\prec e''$ holds.
(Note that the equality does not hold since $e_2\notin L'$.)
Since $st\prec ac=e_2\prec e''$,
\begin{eqnarray*}
e''&=&\max\{e\mid e\in L'-L^*\} \\
   &=&\max\{e\mid e\in (L-ac+st)-L^*\} \\
   &=&\max\{e\mid e\in L-L^*\}.
\end{eqnarray*}
Then $e''$ would have been selected instead of $e_2$ when the parent function $f_1$ is applied to $L$,
which contradicts $e_2=ac$.
Thus, (d) holds.

(``if''-part.)
From (a) and (b), $L=L'-e_1+e_2\in\mathcal{DL}$.
Since $e_1\in L^*$ from (a), (d) implies
\begin{eqnarray*}
e_2&\succ&\max\{e\mid e\in L'-L^*\} \\
    &=& \max\{e\mid e\in (L+e_1-e_2)-L^*\} \\
    &=& \max\{e\mid e\in (L-e_2)-L^*\}.
\end{eqnarray*}
Thus, $e_2=\max\{e\mid e\in L-L^*\}$ holds,
and hence $f_1$ chooses $e_2$ for an edge $ac$ to be deleted from $L$.

From this we have $L-ac=L'-e_1+e_2-ac=L'-e_1$.
Since $e_2\notin L^*$ from (b), (c) implies
\begin{eqnarray*}
e_1&\prec&\min\{e\in L^*-L'\mid L'-e_1+e\in \mathcal{L}\} \\
   &=&\min\{e\in L^*-(L+e_1-e_2)\mid L-ac+e\in \mathcal{L}\} \\
   &=&\min\{e\in L^*-(L+e_1)\mid L-ac+e\in \mathcal{L}\}.
\end{eqnarray*}
Since $e_1\in L^*-L$, $e_1=\min\{e\in L^*-L\mid L-ac+e\in \mathcal{L}\}$.
Thus, $f_1$ chooses $e_1$ for an edge to be added,
and $f_1(L)$ returns $L'$. 
\end{proof}

\begin{proof}[Proof of Lemma \ref{lemma:parent2}]
(``only if''-part.)
Since $f_2(L)=L'$ holds, $e_1$ and $e_2$ must be chosen as $st$ and $ac$ of Definition \ref{def:parent}.
From $e_2=ac$, we have 
\begin{equation}
\label{eq:1dof2}
L-ac=L'-e_1+e_2-ac=L'-e_1.
\end{equation} 
Since $st\in T(L-ac)-(L-ac)$ from Definition \ref{def:parent}
and Fact \ref{fact:under_tri}, $st$ is $F$-legal in $T(L-ac)$ and hence $F$-legal in $T(L-ac+st) (=T(L'))$.
Thus, from $e_1=st$, (a) holds.
Since $e_1$ is $F$-legal in $T(L')$, we have 
\begin{equation}
\label{eq:under_tri}
T(L')=T(L'-e_1)=T(L-ac).
\end{equation}
Let $e'=\min\{e\in T(L')-L'\mid L'-e_1+e\in \mathcal{L}\}$.
Suppose that (c) does not hold.
Then $e'\prec e_1$ holds. 
(Note that the equality does not hold since $e_1\in L'-F$.)
From Eq.(\ref{eq:1dof2}) and (\ref{eq:under_tri}),
\begin{eqnarray*}
e'&=&\min\{e\in T(L')-L'\mid L'-e_1+e\in \mathcal{L}\} \\
  &=&\min\{e\in T(L-ac)-(L-ac+e_1)\mid L-ac+e\in \mathcal{L}\} \\
  &=&\min\{e\in T(L-ac)-(L-ac)\mid L-ac+e\in \mathcal{L}\}.
\end{eqnarray*}
Then, $e'$ would have been selected
when the parent function is applied to $L$, which contradicts $e_1=st$. Hence (c) holds.

Next let us consider $e_2$.
Since $e_2(=ac)$ is $F$-illegal in $T(L)$ from Definition \ref{def:parent},
we have $e_2\notin T(L-ac)=T(L')$.
Thus (b) holds.

Finally, (d) must hold since the parent function remove an lexicographically largest illegal edge in $T(L)$. 

(``if''-part.)
From (a), $e_1$ is $F$-legal in $T(L')$.
Then, we have $T(L')=T(L'-e_1)$.
Suppose that $e_2$ is $F$-legal in $T(L)=T(L'-e_1+e_2)$, 
$T(L'-e_1+e_2)=T(L'-e_1)$ holds and
$e_2$ is in $T(L'-e_1)=T(L')$, which contradicts the condition (b).
Therefore, $e_2$ is $F$-illegal edge in $T(L)$, 
and (d) says that $e_2$ is the lexicographically largest $F$-illegal edge in $T(L)$.
Thus, $f_2$ chooses $e_2$ for an edge $ac$ to be deleted from $L$,
and $L'-e_1=L-ac$ holds.

From $L'-e_1=L-ac$ and the condition (c),
$e_1\prec \min\{e\in T(L-ac)-(L-ac+e_1)\mid L-ac+e\in\mathcal{L}\}$ holds.
Thus, $e_1=\min\{e\in T(L-ac)-(L-ac)\mid L-ac+e \in\mathcal{L}\}$.
(Note that $e_1\in T(L-ac)$, because $T(L-ac)=T(L'-e_1)=T(L')$ and $e_1\in T(L')$ hold from (a).)
Thus, $f_2$ chooses $e_1$ for an edge to be added,
and $f_2(L)$ returns $L'$. 
\end{proof}

By Lemmas \ref{lemma:parent1} and \ref{lemma:parent2},
we will describe an $O(n^3)$ algorithm in the proof of the following Theorem~\ref{theo:time}.
Before it, we give you the simple observation for checking the condition (d) in the Lemma~\ref{lemma:parent2} efficiently:
\begin{observation}
\label{ob:condition_d}
Let $T(L')$ be a $L'$-constrained Delaunay triangulation constrained by edges of a non-crossing Laman framework $L'$,
and let $e_1\in L'$ be a $F$-legal edge in $T(L')$ and $e_2\in K_n-L'$ be an edge intersecting no edge of $L'$.
Then $T(L'-e_1+e_2)=T(L'+e_2)$ holds.
\end{observation}
\begin{proof}
This comes from the following two facts that:
\begin{itemize}
\item[(1)] $T(L')=T(L'-e_1)$ holds, and
\item[(2)] $e_1$ is $F$-legal in $T(L'+e_2)$.
\end{itemize}
The fact (1) clearly hold since $e_1$ is $F$-legal edge in $T(L')$.
Using Lemma~4.2 in \cite{Floriani:1992} the fact (2) immediately follows.
Let us show how to update $T(L')$ to $T(L'+e_2)$.
The edge $e_2$ intersects some edges in $T(L')-L'$, and we denote a set of such edges by $C_{e_2}$.
First, we delete the edges in $C_{e_2}$.
The resulting graph $T(L')-C_{e_2}$ has a hole bounded by the polygon $Q$.
When inserting $e_2$ into $T(L')-C_{e_2}$, $e_2$ splits $Q$ into two polygon $Q_1$ and $Q_2$.
Then, $T(L'+e_2)=(T(L')-C_{e_2})+T(Q_1)+T(Q_2)$ holds from the fact in \cite{Floriani:1992},
i.e. $T(L'+e_2)$ is obtained from $T(L')$ by (Delaunay) triangulating only $Q_1$ and $Q_2$ independently.
This implies that $e_1$ never flips even if $e_1$ is not a constrained edge.
(Note that $e_2$ does not intersect $e_1$ from the lemma assumption.)
Hence, $e_1$ is $F$-legal in $T(L'+e_2)$.  
\end{proof}

\begin{thm}
\label{theo:time}
The set of all $F$-constrained non-crossing Laman frameworks on a given point set can be reported in $O(n^3)$ time 
per output
using $O(n^2)$ space, or $O(n^4)$ time using $O(n)$ space.
\end{thm}
\begin{proof}
As described in Section \ref{sec:main}, we use a linear transformation if necessary to get a unique CDT.
The complexity of testing the uniqueness of a CDT is $O(n^2)$ by simply testing
the circumcircle of each triangle in the CDT
to see there is another point other than vertices of the triangle on the circumcircle.


Given a non-crossing Laman framework $L'\in\mathcal{L}$ and $L'$-constrained Delaunay triangulation $T(L')$.
The algorithm will check $f_1(Adj(L',e_1,e_2))=L'$ or $f_2(Adj(L',e_1,e_2))=L'$ at line~10 depending on the edge pair $(e_1,e_2)$.
Here we will show that each condition in Lemmas \ref{lemma:parent1} and \ref{lemma:parent2} can be checked in $O(1)$ time
for each of the $O(n^3)$ edge pairs $(e_1,e_2)$ by the following way.

First, for all edges $e_2\in {\sf elist}_{K_n}$,
we calculate the number of edges $e_1\in L'$ intersecting $e_2$, which we denote by {\sf cross\_n}$(e_2, L')$.
If {\sf cross\_n}$(e_2, L')>1$, we delete $e_2$ from ${\sf elist}_{K_n}$ since $L'-e_1+e_2$ is never non-crossing for any $e_1\in {\sf elist}_{L'}$.
If {\sf cross\_n}$(e_2, L')=1$, we associate a pointer of the edge $e_1$ intersecting $e_2$ with $e_2$, 
and we denote such $e_1$ by {\sf cross\_e}$(e_2, L')$. 

Next, for each $e_1\in {\sf elist}_{L'}$, we attach two flags to $e_1$ which represent that
$e_1$ satisfies the condition (a) of Lemmas \ref{lemma:parent1} and \ref{lemma:parent2}, respectively.
This preprocessing can be done in $O(n)$ time by simply scanning all edges of ${\sf elist}_{L'}$.
By this we can check the condition (a) in Lemmas \ref{lemma:parent1} and \ref{lemma:parent2} in $O(1)$ time.
Similarly, we attach two flags to $e_2\in {\sf elist}_{K_n}$ which represent that $e_2$ satisfies the condition (b) of 
Lemmas \ref{lemma:parent1} and \ref{lemma:parent2}.
This process can be done in $O(n^2)$ time.
By this we can check the condition (b) in Lemmas \ref{lemma:parent1} and \ref{lemma:parent2} in $O(1)$ time.
We can calculate the lexicographically largest edge in $L'-L^*$ in $O(n)$ time.
By this we can check the condition (d) in Lemmas \ref{lemma:parent1} in $O(1)$ time.

Now let us consider how to identify a set of edges $e_2\in {\sf elist}_{L'}$ satisfying the condition (d) 
in Lemma~\ref{lemma:parent2} in $O(n^3)$ time with $O(n^2)$ space.
(In the case of $O(n^4)$ time algorithm this process must be skipped, and the condition (d) in Lemma~\ref{lemma:parent2} will be checked
simply by updating $T(L')$ to $T(L'-e_1+e_2)$ for each pair $(e_1,e_2)$ using $O(n)$ time and $O(n)$ space 
by applying the algorithm by Chin and Wang\cite{chin:wang:1998}.)
It can be done regardless of the removing edge $e_1$ if the condition (a) in Lemma~\ref{lemma:parent2} holds.
From Observation~\ref{ob:condition_d}, we can say that the condition~(d) holds if and only if 
$e_2$ is lexicographically largest $F$-illegal edge in $T(L+e_2)$ when {\sf cross\_n}$(e_2, L')=0$.
And it is sufficient to check the condition~(d) only in $T(L-${\sf cross\_e}$(e_2, L')+e_2)$ when {\sf cross\_n}$(e_2, L')=1$.
Updating the Delaunay triangulation takes $O(n)$ time
(see \cite{Anglada:1997,Floriani:1992,chin:wang:1998} for a linear time update of the constrained Delaunay triangulation).
Thus we can attach a flag to $e_2\in {\sf elist}_{K_n}$ in $O(n)$ time 
which represents whether $e_2$ satisfies the condition (d) of Lemma~\ref{lemma:parent2} or not,
thus this preprocessing for all edges in ${\sf elist}_{K_n}$ takes $O(n^3)$ time.

By using the above mentioned data, we will show that for a fixed $e_1\in {\sf elist}_{L'}$,
the inner while-loop from line 6 to 16 can be executed in $O(n^2)$.
In order to efficiently test the condition (c) of Lemmas \ref{lemma:parent1} and \ref{lemma:parent2}, we prepare   
the data structure proposed by  Lee, Streinu and Theran \cite{lee:streinu:2005,lee:streinu:theran:2005} in $O(n^2)$ time
for maintaining  {\it rigid components} obtained by deleting $e_1$.
This data structure supports a {\it pair-find} query which determines whether two vertices are spanned by a common component 
in $O(1)$ time using $O(n^2)$ space, or $O(n)$ time with $O(n)$ space.
From this, we can calculate $Adj(L',e_1,e_2)$ (i.e., determine whether $L'-e_1+e_2\in {\cal L}$) 
in $O(1)$ time with $O(n^2)$ space, or $O(n)$ time with $O(n)$ space, for each edge $e_2\in {\sf elist}_{K_n}$. 
Also, we can compute $e'=\min\{e\in L^*-L'\mid L'-e_1+e\in\mathcal{L}\}$ and 
$e''=\min\{e\in T(L')-L'\mid L'-e_1+e \in\mathcal{L}\}$ in $O(n)$ time with $O(n^2)$ space, or $O(n^2)$ time with $O(n)$ space.
From $e'$ and $e''$ we can check condition (c) in Lemmas \ref{lemma:parent1} and \ref{lemma:parent2} in $O(1)$ time.


We have to determine  which of $f_1(Adj(L',e_1,e_2))=L'$ and 
$f_2(Adj(L',e_1,e_2))=L'$ should be checked for every  edge 
$e_2\in {\sf elist}_{K_n}$.
This is simply done by checking  whether $e_2\in T(L^*)-L^{*}$ or $e_2\in K_n-T(L')$.
When $e_2\in T(L^*)-L^*$, we check whether $f_1(Adj(L',e_1,e_2))=L'$ holds or not.
When $e_2\in K_n-T(L')$, we check whether $f_2(Adj(L',e_1,e_2))=L'$ holds or not.
Both can be done in $O(1)$ time with $On^2$ space, or $O(n)$ time with $O(n)$ space.

By using the above mentioned data structure for maintaining the rigid components,
we can perform both parent function and adjacency function in $O(n^2)$ time with $O(n^2)$ space, or $O(n^3)$ time with $O(n)$ space.
Thus, we have an $O(n^3)$ algorithm using $O(n^2)$ space, or $O(n^4)$ algorithm using $O(n)$ space.
\end{proof}

\begin{figure}[h]
\centering
\includegraphics[width=0.97\textwidth]{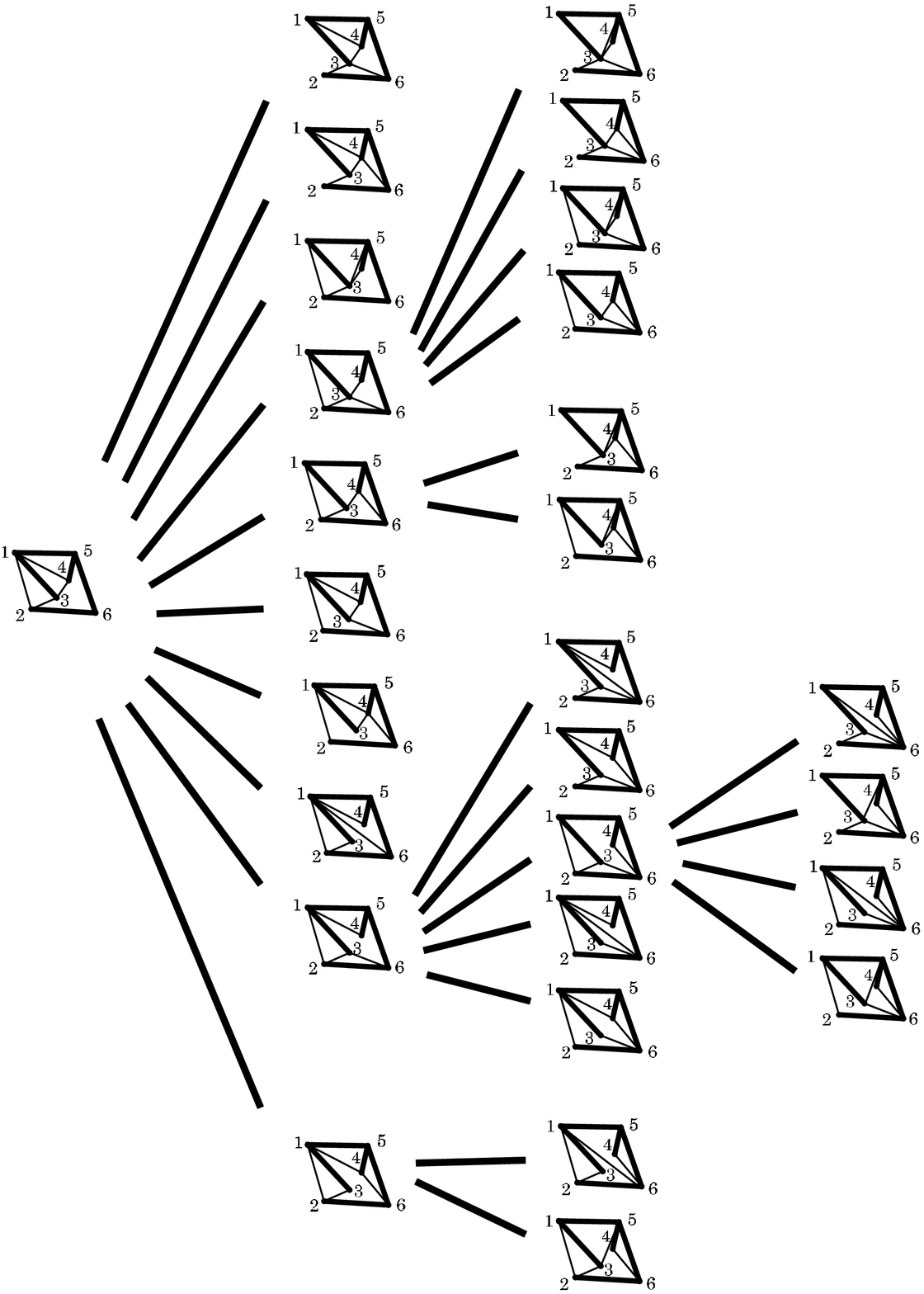}
\caption{An example of the search tree of our algorithm on six points, where $F=\{13,15,26,45,56\}$.}
\label{fig:search_tree}
\end{figure}

Figure \ref{fig:search_tree} illustrates an example of the search tree of a set of $F$-constrained non-crossing Laman frameworks on six points, 
where $F=\{13,15,26,45,56\}$.

\section{ Conclusions}

We have presented an algorithm for enumerating all the constrained non-crossing Laman frameworks. 
%
We note in passing that the techniques in this paper can also be used to 
generate all $F$-constrained non-crossing spanning trees of a point set
since they also form bases of the graphic matroid on any triangulation of $P$.
The unconstrained case was considered in \cite{avis:fukuda:reverse-search:1996}.

\small
\bibliographystyle{abbrv}

\end{document}